\documentclass[11pt]{amsart} 

\usepackage[lmargin=1in,rmargin=1in,tmargin=1in,bmargin=1in]{geometry}
\usepackage[ps,all,arc,rotate]{xy}
\usepackage{graphicx, float, epstopdf}
\usepackage{hyperref, color}
\usepackage{centernot}
\usepackage{fancyhdr}
\usepackage[utf8]{inputenc}
\usepackage{amsfonts,amssymb,amsmath,amsthm,mathrsfs}
\usepackage{graphics, setspace}
\usepackage{braket}
\usepackage{mathtools}
\usepackage{tikz}

\usepackage[
backend=biber,
style=alphabetic,
]{biblatex}
\numberwithin{equation}{section}
\numberwithin{figure}{section}
  
\allowdisplaybreaks[4]         

\newtheorem{lemma}{Lemma}[section]
\newtheorem{theorem}{Theorem}[section]

\newtheorem{corollary}{Corollary}[section]
\theoremstyle{definition}

\newtheorem{remark}{Remark}[section]

\newcommand{\bZ}{\mathbb{Z}}
\newcommand{\bQ}{\mathbb{Q}}
\newcommand{\bA}{\mathbb{A}}
\newcommand{\Sym}{\text{Sym}}
\newcommand{\cF}{\mathcal{F}}
\newcommand{\ve}{\varepsilon}

\newcommand{\bR}{\mathbb{R}}
\newcommand{\bC}{\mathbb{C}}

\begin{document}
\title{On Siegel zeros of symmetric power L-functions} 
\author{Shifan Zhao}
\address{Department of Mathematics, The Ohio State University, United States}
\email{zhao.3326@osu.edu}

\maketitle
	
\begin{abstract}
Let $f$ be a holomorphic cusp form of even weight $k$ for the modular group SL$(2,\bZ)$, which is assumed to be a common eigenfunction for all Hecke operators. For positive integer $n$, let $\Sym^n(f)$ be the symmetric $n^\text{th}$ power lifting of $f$, which was shown by Newton and Thorne \cite{NT1} to be automorphic and cuspidal. In this paper, we construct certain auxiliary $L$-functions to show that Siegel zeros of $L(s,\Sym^n(f))$ do not exist, for each given $n$, utilizing the above functorality result. As an application, we give a lower bound of those symmetric power $L$-functions at $s=1$ of logarithm power type.
\end{abstract}
	
\section{Introduction} \label{sec:introduction}

Siegel zeros of $L$-functions are real zeros that are very close to 1. Precisely, let $L(s,f)$ be an $L$-function (In this paper, by $L$-functions we mostly refer to either automorphic $L$-functions or Rankin-Selberg type $L$-functions) and let $c>0$ be a positive number. A Siegel zero/Landau-Siegel zero/exceptional zero of $L(s,f)$ relative to $c$ is a real zero of $L(s,f)$ that lies in the interval 
$$(1-\frac{c}{\log(Q_f)},1)$$
where $Q_f$ is the analytic conductor of $L(s,f)$, introduced by Iwaniec and Sarnak in \cite{IS}. In practice we will consider families of $L$-functions and we will be interested in finding an absolute and effective constant $c$ that works for every member in a family. We thus introduce the following terminology: Let $\cF$ be a family of $L$-functions. We say $\cF$ admits no Siegel zeros, if there exists a positive effective constant $c_\cF>0$ that depends only on $\cF$, such that each member in $\cF$ has no Siegel zeros relative to $c_\cF$. 
\par
It is conjectured that Siegel zeros of automorphic $L$-functions should not exist, for such zeros, if exist, will be close to 1 and thus not on the critical line, which would contradict to the Generalized Riemann Hypothesis (GRH).
\par
It is an important topic in analytic number theory to rule out the hypothetical exceptional zeros for various families of $L$-functions, not only because they are potential counterexamples of GRH, but also because of their intrinsic significance. The problem of Siegel zeros is closely connected to the class number problem, the distribution of primes in arithmetic progressions, and the problem of estimating Fourier coefficients of classical modular forms, to name a few. See \cite{D} and \cite{HL} for a discussion of these connections.
\par
The problem of Siegel zeros of Dirichlet $L$-functions remains open today and is considered as one of the deepest and most difficult problems in number theory. In 1994, Goldfeld, Hoffstein and Lieman  \cite{GHL} proved that $L$-functions of Gelbart-Jacquet symmetric square lifts of non-dihedral GL(2) cusp forms admit no Siegel zeros. Later in 1995 Hoffstein and Ramakrishnan \cite{HR} showed that $L$-functions of GL(2) cusp forms over any number field $F$ admit no Siegel zeros. As a corollary, in the case when $F = \bQ$, they gave a Siegel-type lower bound to special $L$-values $L(1,\pi)$, where $\pi$ is a GL(2) cusp form. They also studied $L$-functions of GL(3) cusp forms, and conditionally eliminated their Siegel zeros, under certain hypothesis on twisted symmetric square $L$-functions. That hypothesis was later proved by Banks \cite{B} in 1997, and the non-existence of Siegel zeros on GL(3) was established unconditionally. Siegel zeros of Rankin-Selberg type $L$-functions were studied by Ramakrishnan and Wang \cite{RW} in 2003. They proved that Rankin-Selberg type $L$-functions $L(s,\pi \times \pi^\prime)$ admit no Siegel zeros, where $\pi$ and $\pi^\prime$ are GL(2) cuspidal representation over any number field, with a few exception cases. They also eliminated Siegel zeros of symmetric fourth power $L$-functions of a self-dual GL(2) cusp form $\pi$, not of solvable polyhedral type. This gives non-existence of Siegel zeros of $L$-functions of the form $L(s,\Sym^2(\pi) \times \Sym^2(\pi))$, as a direct corollary. Siegel zeros on the Galois side were studied by Wang in 2003 \cite{W}. There he proved that for a cusp form $\pi$ on GL(2) over a number field $F$ of strongly icosahedral type (that is, it is attached to a strongly modular icosahedral Galois representation), and an idele character $\chi$, the twisted symmetric power $L$-functions $L(s,\Sym^m(\pi) \otimes \chi)$ admit no Siegel zeros, unless in the reduced case that $\Sym^m(\pi) \otimes \chi$ has a constituent of a trivial or quadratic character. 
\par
In this paper we study symmetric power $L$-functions of holomorphic Hecke eigenforms of SL$(2,\bZ)$, the automorphy of which was established by Newton and Thorne \cite{NT1} in 2021. We prove their Siegel zeros do not exist. As a consequence, we give a sharp lower bound of these symmetric power $L$-functions at $s=1$. Our main results are summarized below: 

\begin{theorem}
Let $f$ be a normalized holomorphic Hecke eigenform of SL$(2,\bZ)$ of weight $k \geq 2$. Let $n \geq 1$ be a positive integer. Then the symmetric $n^{\text{th}}$ power $L$-functions $L(s,\Sym^n(f))$ admit no Siegel zeros. Explicitly, there exists an effective constant $c_n>0$, depending only on $n$, such that $L(s,\Sym^n(f))$ have no real zeros in the interval
$$(1-\frac{c_n}{\log k},1)$$
\end{theorem}

\begin{remark}
In fact, the dependence of $c_n$ on $n$ can be made explicit. We will prove in section 3 that there exists an absolute effective constant $c>0$ such that $L(s,\Sym^n(f))$ has no real zeros in $(1-\frac{c}{n^4\log(nk)},1)$.
\end{remark}

\begin{remark}
In this paper, we work with level 1 holomorphic cusp forms for simplicity. In fact, in the second part of Newton and Thorne's paper \cite{NT2}, symmetric power liftings of all non-CM regular algebraic cuspidal automorphic representations $\pi$ of GL$_2(\bA_\bQ)$ were established. We expect similar results for those symmetric power $L$-functions. We will study this issue carefully in upcoming papers.
\end{remark}

As a corollary, we give the following lower bound for the special $L$-values $L(1,\Sym^n(f))$:

\begin{corollary}
Let $f$ be a normalized holomorphic Hecke eigenform of SL$(2,\bZ)$ of weight $k \geq 2$. Let $n \geq 1$ be a positive integer. Then for any $\ve>0$ the following lower bound holds:
$$L(1,\Sym^n(f)) \gg_{n,\ve} \frac{1}{\log(k)^{2n+2+\ve}}$$
where the implied constant is effective and depends only on $n$ and $\ve$.
\end{corollary}

\begin{remark}
Before we state the proof in section 4, a few remarks about Corollary 1.1 should be made here. First, the key ingredient in the proof is to take an appropriate root of the symmetric power $L$-function inside its zero-free region. The idea of taking square root was first came up with by W. Luo, as indicated in the remark of Theorem C in \cite{HR}. Using similar approach, Cogdell and Michel \cite{CM} obtained parallel results for holomorphic newforms of weight 2 and square-free level $q$, in the $q$-aspect, under the hypothesis that symmetric power $L$-functions of those forms are automorphic and Siegel zeros do not exist. For completeness we state a proof here.
\end{remark}

\section{Preliminaries}

In this section we review some facts that will be essential to our proof, as well as fix some notations.
\par
Let $\pi$ be a cuspidal automorphic representation of GL$_d(\bA_\bQ)$ with unitary central character. There is an $L$-function $L(s,\pi)$ attached to it, which admits a Dirichlet series representation 
$$L(s,\pi) = \sum_{n=1}^\infty \frac{\lambda_\pi(n)}{n^s}, \sigma > 1$$
and an Euler product expansion
$$L(s,\pi) = \prod_{p} \prod_{i=1}^d (1-\frac{\alpha_{\pi,i}(p)}{p^s})^{-1}, \sigma > 1$$
where both the series and the infinite product converge absolutely. 
The $L$-function $L(s,\pi)$ extends to an entire function on the whole complex plane, satisfying a functional equation of the following type:
$$q_\pi^{\frac{s}{2}}L_\infty(s,\pi)L(s,\pi) = \ve_\pi q_\pi^{\frac{1-s}{2}}L_\infty(1-s,\tilde{\pi})L(1-s,\tilde{\pi})$$
where $q_\pi \geq 1$ is the arithmetic conductor of $\pi$, $\ve_\pi$ has norm 1, and $\Tilde{\pi}$ is the contragredient of $\pi$. Here $L_\infty(s,\pi)$ is a product of gamma factors:
$$L_\infty(s,\pi) = \prod_{j=1}^d \Gamma_\bR(s-\mu_{\pi,j})$$
where $\Gamma_\bR = \pi^{-\frac{s}{2}}\Gamma(\frac{s}{2})$.
\par
For finitely many cuspidal representations $\pi_1, \pi_2, \dots, \pi_r$ of degree $d_1, d_2, \dots, d_r$ respectively, Langlands found a distinguished automorphic representation $\pi$ of degree $d = \sum_{i=1}^r d_i$, whose $L$-function is the product of the cuspidal ones:
$$L(s,\pi) = \prod_{i=1}^r L(s,\pi_i)$$
This representation $\pi$ is called the isobaric sum of $\pi_1, \dots, \pi_r$, and is denoted by $\pi_1 \boxplus \cdots \boxplus \pi_r$. Moreover, an automorphic representation that is isomorphic to an isobaric sum of finitely many cuspidal representations is called an isobaric representation. 
\par
Let $\pi_1$, $\pi_2$ be two cuspidal representations of degree $d_1,d_2$ respectively. There is an $L$-function $L(s,\pi_1 \times \pi_2)$ attached to them, called the Rankin-Selberg convolution, with a Dirichlet series expansion
$$L(s,\pi_1 \times \pi_2) = \sum_{n=1}^\infty \frac{\lambda_{\pi_1 \times \pi_2}(n)}{n^s}, \sigma > 1$$
and an Euler product
$$L(s,\pi_1 \times \pi_2) = \prod_p L_p(s,\pi_1 \times \pi_2)$$
where for $p$ not dividing the conductor of $\pi_1$ and $\pi_2$ (such a prime is called ramified), the local $L$-function $L_p(s,\pi_1 \times \pi_2)$ takes the following form:
$$L_p(s,\pi_1 \times \pi_2) = \prod_{i=1}^{d_1}\prod_{j=1}^{d_2} (1-\frac{\alpha_{\pi_1,i}(p)\alpha_{\pi_2,j}(p)}{p^s})^{-1}$$
Moreover, $L(s,\pi_1 \times \pi_2)$ extends to $\bC$, admitting a simple pole at $s=1$ precisely when $\pi_1 \cong \overline{\pi_2}$. And it also satisfies a functional equation.
\par
Rankin-Selberg convolutions can be defined for isobaric representations, by formally distributing $\times$ with $\boxplus$.
\par
The following lemma about Rankin-Selberg $L$-functions of an isobaric representation with its contragredient was proved in \cite{HR}

\begin{lemma}
Let $\Pi$ be an isobaric representation of GL$_d(\bA_\bQ)$. Then the Rankin-Selberg $L$-function $L(s,\Pi \times \Tilde{\Pi})$ has non-negative Dirichlet series coefficients. 
\end{lemma}

In the proof we will deal with $L$-functions with non-negative Dirichlet series coefficients, whose Siegel zeros were studied in \cite{GHL}. We will use a slightly more refined version in \cite{IK}.

\begin{lemma}
Let $L(s,f)$ be an $L$-function of degree $d$ with non-negative Dirichlet series coefficients. Suppose $L(s,f)$ has a pole at $s=1$ of order $r \geq 1$. Then there exists an absolute and effective constant $c>0$ such that $L(s,f)$ has at most $r$ real zeros (counting multiplicity) in the interval
$$(1-\frac{c}{d(r+1)\log Q_f},1)$$
where $Q_f$ is the analytic conductor of $L(s,f)$.
\end{lemma}

Lastly, we recall the definition of symmetric power $L$-functions. Let $f$ be a holomorphic Hecke eigenform of SL$(2,\bZ)$ of weight $k$, with Fourier expansion
$$f(z) = \sum_{m=1}^\infty \lambda_f(m)m^{\frac{k-1}{2}}e(mz)$$
and Hecke $L$-function
$$L(s,f) = \sum_{m=1}^\infty \frac{\lambda_f(m)}{m^s} = \prod_p(1-\frac{\alpha_p}{p^s})^{-1}(1-\frac{\beta_p}{p^s})^{-1}, \sigma > 1$$
For $n \geq 0$ the symmetric $n^{\text{th}}$ power $L$-function of $f$ is defined by the Euler product for $\sigma > 1$:
$$L(s,\Sym^n(f)) = \prod_p\prod_{j=0}^n(1-\frac{\alpha_p^j\beta_p^{n-j}}{p^s})^{-1}$$
which was proved in \cite{NT1} to be cuspidal automorphic. So in particular they extend to entire functions, unless $n=0$, where $L(s,\Sym^0(f)) = \zeta(s)$ is the Riemann zeta function.

\section{Proof of Theorem 1.1}
We begin with the following decomposition lemma, which states that Rankin-Selberg $L$-function of two symmetric powers decomposes into a product of symmetric power $L$-functions.

\begin{lemma}
Let $f$ be as in Theorem 1.1. Let $n \geq 1$ and $r \geq 0$ be integers. Then
$$L(s,\Sym^n(f) \times \Sym^{n+r}(f)) = \prod_{i=0}^n L(s,\Sym^{2i+r}(f))$$
\end{lemma}

\begin{remark}
This lemma is probably known to experts but we cannot find any literature that clearly states it. So we give our proof here. 
\end{remark}

\begin{proof}
We verify this identity by showing the local parameters of both sides are identical. Let $p$ be a prime and let $\alpha_p$, $\beta_p$ be the local parameters of $f$ at $p$. Note that $\alpha_p\beta_p = 1$. The (multi)set of local parameters at $p$ of the left side is
$$\{\alpha_p^i\beta_p^{n-i}\alpha_p^{j}\beta_p^{n+r-j}: 0 \leq i \leq n, 0 \leq j \leq n+r\}$$
In view of the relation $\alpha_p\beta_p = 1$, the multiset above can be written as
$$\{\beta_p^{2n+r-2i-2j}:0 \leq i \leq n, 0 \leq j \leq n+r\}$$
The multiset of local parameters at $p$ of the right side is
$$\{\alpha_p^{l}\beta^{2m+r-l}:0 \leq m \leq n, 0 \leq l \leq 2m+r\}$$
which can be written as
$$\{\beta_p^{2m+r-2l}:0 \leq m \leq n, 0 \leq l \leq 2m+r\}$$
Thus is suffices to show the following equality
$$\{2n+r-2i-2j:0 \leq i \leq n, 0 \leq j \leq n+r\} = \{2m+r-2l:0 \leq m \leq n, 0 \leq l \leq 2m+r\}$$
By shifting by $r$ and then divided by 2, it suffices to show the two multisets $A_{n,r}$ and $B_{n,r}$ are equal:
\begin{align*}
A_{n,r} &= \{n-i-j:0 \leq i \leq n, 0 \leq j \leq n+r\} \\
B_{n,r} &= \{m-l:0 \leq m \leq n, 0 \leq l \leq 2m+r\}
\end{align*}
We fix $r \geq 0$ and prove by induction on $n \geq 1$.
\begin{itemize}
\item
The base case is $n=1$, where it is easy to verify that $A_{1,r}$ and $B_{1,r}$ are both equal to 
$$\{1,0^2,(-1)^2,\dots,(-r)^2, 1-r\}$$
where $x^2$ indicates that the element $x$ appears twice.
\item
Assume $A_{n,r} = B_{n,r}$ for some $n \geq 1$. We need to prove $A_{n+1,r} = B_{n+1,r}$. For $A_{n+1,r}$ we have
\begin{align*}
A_{n+1,r} &= \{n+1-i-j:0 \leq i \leq n+1, 0 \leq j \leq n+1+r\} \\
&= (A_{n,r}+1) \cup \{-i-r:0 \leq i \leq n\} \cup \{-j:0 \leq j \leq n+r+1\} \\
&= (A_{n,r}+1) \cup \{-n-r-1,(-n-r)^2,\dots,(-r)^2,-r+1,\dots,0\}
\end{align*}
For $B_{n+1,r}$ we have
\begin{align*}
B_{n+1,r} &= \{m-r:0 \leq m \leq n+1, 0 \leq l \leq 2m+r\} \\
&= \{m-l:1 \leq m \leq n+1, 0 \leq l \leq 2m+r\} \cup \{-l:0 \leq l \leq r\} \\
&= \{m-l+1:0 \leq m \leq n,0 \leq l \leq 2m+r+2\} \cup \{-l:0 \leq l \leq r\} \\
&= (B_{n,r}+1) \cup \{-m-r,-m-r-1:0 \leq m \leq n\} \cup \{-l:0 \leq l \leq r\} \\
&= (B_{n,r}+1) \cup \{-n-r-1,(-n-r)^2,\dots,(-r)^2,-r+1,\dots,0\}
\end{align*}
where in the third line we changed $m$ to $m+1$. By induction hypothesis we have $A_{n+1,r} = B_{n+1,r}$, and the proof is finished.
\end{itemize}
\end{proof}

\begin{proof}
For the proof of Theorem 1.1, we need to construct certain auxiliary $L$-functions. The construction depends on the parity of $n$. 
\par
First assume $n$ is even. Let $\Pi = 1 \boxplus \Sym^n(f)$, which is isobaric and self-dual. By Lemma 2.1 the auxiliary $L$-function $D(s) = L(s,\Pi \times \Pi)$ has non-negative Dirichlet series coefficients. The function $D(s)$ is of degree $(n+2)^2$ and has a double pole at $s=1$. Therefore by Lemma 2.2 there exists absolute and effective constant $c>0$
such that $D(s)$ has at most two real zeros in the interval
$$I = (1-\frac{c}{3(n+2)^2\log Q_{\Pi \times \Pi}},1)$$
Now by our lemma 3.1, $D(s)$ can be decomposed as follows:
\begin{align*}
D(s) &= \zeta(s)L(s,\Sym^n(f))^2L(s,\Sym^n(f) \times \Sym^n(f)) \\
&= \zeta(s)L(s,\Sym^n(f))^2\prod_{i=0}^nL(s,\Sym^{2i}(f)) \\
&= \zeta(s)^2L(s,\Sym^n(f))^3\prod_{1 \leq i \leq n, i \neq \frac{n}{2}}L(s,\Sym^{2i}(f))
\end{align*}
Since $\zeta(s)$ has only one simple pole at $s=1$ and each symmetric power $L$-function $L(s,\Sym^{2i}(f))$ ($1 \leq i \leq n, i \neq \frac{n}{2}$) is entire, any zero of $L(s,\Sym^n(f))$ in the interval $I$ would be a zero of $D(s)$ of order at least 3. Since $D(s)$ has at most two zeros in $I$, $L(s,\Sym^n(f))$ has no zeros in $I$ at all.
\par
We now estimate logarithm of the analytic conductor $\log (Q_{\Pi \times \Pi})$, in terms of $n$ and $k$. Since $L(s,\Pi \times \Pi)$ can be written as a product of symmetric power $L$-functions, it suffices to compute the analytic conductor of a symmetric power $L$-function. The gamma factors $L_\infty(s,\Sym^n(f))$ have been worked out for example in \cite{LW} and \cite{CM}. They takes different forms, depending on the parity of $n$. Explicitly, when $n=2m+1$ is odd, we have
$$L_\infty(s,\Sym^n(f)) = \prod_{j=0}^m \Gamma_\bC(s+(j+\frac{1}{2})(k-1))$$
And when $n=2m$ is even, we have
$$L_\infty(s,\Sym^n(f)) = \Gamma_\bR(s+\delta_{2 \nmid m}) \prod_{j=1}^m \Gamma_\bC(s+j(k-1))$$
where $\Gamma_\bR(s) = \pi^{-\frac{s}{2}}\Gamma(\frac{s}{2})$, $\Gamma_\bC = 2(2\pi)^{-s}\Gamma(s)$, and $\delta_{2 \nmid m} = 1$ if $2 \nmid m$, and $0$ otherwise. Thus when $n=2m$ is even, the following estimate holds:
\begin{align*}
Q_{\Sym^n(f)} &= (1+\delta_{2 \nmid m})\prod_{j=1}^m (j(k-1)+1)(j(k-1)+2) \\
&\ll \prod_{j=1}^m j^2(k-1)^2 \\
&\ll (k-1)^{2m}(m!)^2 \\
&= (k-1)^n(m!)^2
\end{align*} 
Taking logarithm, we have
\begin{align*}
\log(Q_{\Sym^n(f)}) &\ll n\log(k-1)+2\log(m!) \\
&\ll n\log k + 2m\log m \\
&\ll n\log(\frac{nk}{2})
\end{align*}
using Stirling's approximation. For $n$ odd we run a similar argument and get the same estimate. 
\par
Now note that $L(s,\Pi \times \Pi)$ is decomposed into a product of $(n+4)$ symmetric power $L$-functions, among which the highest power is $2n$. So we have
\begin{align*}
\log(Q_{\Pi \times \Pi}) &\ll 2(n+4)n\log(nk) \\
&\ll n^2\log(nk)
\end{align*}
Therefore, $L(s,\Sym^n(f))$ does not vanish in the interval
$$(1-\frac{c}{3(n+2)^2n^2\log(nk)},1)$$
By changing $c$ if necessary, we proved that $L(s,\Sym^n(f))$ has no zeros in the interval
$$(1-\frac{c}{n^4\log(nk)},1)$$
as claimed, and all constants involved are absolute and effective.
\par
In the other case where $n$ is odd, we consider the isobaric representation $\Pi = 1 \boxplus \Sym^n(f) \boxplus \Sym^{n+1}(f)$ instead, and still examine the auxiliary $L$-function $D(s) = L(s,\Pi \times \Pi)$. The following decomposition of $D(s)$ holds:
\begin{align*}
D(s) &= \zeta(s)L(s,\Sym^n(f))^2L(s,\Sym^{n+1}(f))^2L(s,\Sym^n(f)\times \Sym^n(f)) \\
&\times L(s,\Sym^n(f) \times \Sym^{n+1}(f))^2L(s,\Sym^{n+1}(f) \times \Sym^{n+1}(f)) \\
&= \zeta(s)L(s,\Sym^n(f))^2L(s,\Sym^{n+1}(f))^2\prod_{i=0}^nL(s,\Sym^{2i}(f)) \\
&\times \prod_{j=0}^nL(s,\Sym^{2j+1}(f))^2\prod_{k=0}^{n+1}L(s,\Sym^{2k}(f)) \\
&= \zeta(s)^3L(s,\Sym^n(f))^4L(s,\Sym^{n+1}(f))^2\prod_{i=1}^nL(s,\Sym^{2i}(f)) \\
&\times \prod_{0 \leq j \leq n, j \neq \frac{n-1}{2}}L(s,\Sym^{2j+1}(f))^2\prod_{k=1}^{n+1}L(s,\Sym^{2k}(f))
\end{align*}
which has a pole of order 3 at $s=1$ and is divided by $L(s,\Sym^n(f))$ with multiplicity 4, with other terms having at most one pole at $s=1$.  The rest of argument is similar to the even case (using the fact that $4>3$)
\end{proof}

\begin{remark}
The general principle that functorality would imply non-existence of Siegel zero was first pointed out in \cite{HR}, section 4. We worked out the details here in the special case of symmetric power $L$-functions. Moreover, we give an explicit (and effective) zero-free interval for those $L$-functions, not only in terms of weight $k$, but also in terms of power $n$.
\end{remark}

\section{Proof of Corollary 1.1}

\begin{proof}
Using Theorem 1.1, combined with the standard zero-free region of automorphic $L$-functions (see, for example, Theorem 5.10 in \cite{IK}), we have the following zero-free region for $L(s,\Sym^n(f))$:
$$R = \{s:\sigma > 1-\frac{2c_n}{\log(k(|t|+1))}\}$$
where $c_n>0$ depends on $n$. In this region we may take the $(n+1)^\text{th}$-root of $L(s,\Sym^n(f))$, which is holomorphic. We then consider the auxiliary $L$-function 
$$D(s) = \zeta(s)L(s,\Sym^n(f))^{\frac{1}{n+1}}$$
which has a simple pole at $s=1$ with residue $L(s,\Sym^n(f))^{\frac{1}{n+1}}$. We compute the logarithm of $D(s)$ for $\sigma>1$:
\begin{align*}
\log D(s) &= \log \zeta(s) + \frac{1}{n+1} \log L(s,\Sym^n(f)) \\
&= \sum_p \log (1-\frac{1}{p^s})^{-1} + \frac{1}{n+1}\sum_p\sum_{j=0}^n \log (1-\frac{\alpha_p^j \beta_p^{n-j}}{p^s})^{-1} \\
&= \sum_p\sum_{m=1}^\infty \frac{1}{p^{ms}} + \frac{1}{n+1}\sum_p\sum_{j=0}^n\sum_{m=1}^\infty \frac{\alpha_p^{mj}\beta_p^{m(n-j)}}{p^{ms}} \\
&= \frac{1}{n+1}\sum_p\sum_{m=1}^\infty \frac{(n+1)+\sum_{j=0}^n\alpha_p^{mj}\beta_p^{m(n-j)}}{p^{ms}}
\end{align*}
where for each $m \geq 1$, the coefficient
$$(n+1)+\sum_{j=0}^n\alpha_p^{mj}\beta_p^{m(n-j)}$$
is a symmetric polynomial in $\alpha_p$ and $\beta_p$. Hence it can be expressed as a polynomial of two elementary symmetric polynomials $\alpha_p+\beta_p = \lambda_f(p)$ and $\alpha_p\beta_p = 1$, where both of them are real, due to Hecke theory. So the coefficient is also real. Moreover, by Deligne's bound we have 
$$|\sum_{j=0}^n\alpha_p^{mj}\beta_p^{m(n-j)}| \leq \sum_{j=0}^n |\alpha_p|^{mj}|\beta_p|^{m(n-j)} = n+1$$
So the coefficient $(n+1)+\sum_{j=0}^n\alpha_p^{mj}\beta_p^{m(n-j)}$ is non-negative. Therefore $\log D(s)$ is a Dirichlet series with non-negative coefficients, and so is $D(s)$. We write
$$D(s) = \sum_{m=1}^\infty \frac{\lambda_D(m)}{m^s}, \lambda_D(1) = 1, 
\lambda_D(m) \geq 0$$
for $\sigma > 1$. Let $x>2$ and $\frac{1}{2} < \beta < 1$ and consider the following integral:
$$I = \frac{1}{2\pi i}\int_{(2)} \frac{D(s+\beta)x^s}{s(s+1)(s+2)}ds$$
where $(2)$ is the vertical line with $\sigma = 2$, pointing upwards. We can expand the series to get
\begin{align*}
I &= \frac{1}{2\pi i}\int_{(2)} \sum_{m=1}^\infty \frac{\lambda_D(m)}{m^{s+\beta}}\frac{x^s}{s(s+1)(s+2)}ds \\
&= \sum_{m=1}^\infty \frac{\lambda_D(m)}{m^\beta} \frac{1}{2\pi i}\int_{(2)} \frac{(x/m)^s}{s(s+1)(s+2)}ds \\
&= \sum_{m<x}\frac{\lambda_D(m)}{m^\beta}\frac{1}{2}(1-\frac{m}{x})^2 \\
&\geq \frac{1}{2}(1-\frac{1}{x})^2 \\
&\gg 1
\end{align*}
Here we used the following integral formula (see \cite{HL}):
$$\frac{1}{2\pi i}\int_{(2)} \frac{x^s}{s(s+1) \cdots (s+r)}ds = \begin{cases}
\frac{1}{r!}(1-\frac{1}{x})^r & x>1 \\
0 & 0<x<1 \\
\end{cases}$$
On the other hand, we may switch the Contour to $\gamma$ and pick up residues, where $\gamma = \gamma_1 \cup \gamma_2 \cup \gamma_3$ is defined as follows:
\begin{align*}
    \gamma_1 &= \{s: \sigma = 1-\beta, |t| \geq k-1\} \\
    \gamma_2 &= \{s: \frac{c_n}{(\log k)^{2+\epsilon}} - \frac{c_n}{\log k} \leq \sigma \leq 1-\beta, |t| = k-1\} \\
    \gamma_3 &= \{s: \sigma = \frac{c_n}{(\log k)^{2+\epsilon}} - \frac{c_n}{\log k}, |t| \leq k-1\}
\end{align*}
We choose $\beta = 1-\frac{c_n}{\log(k)^{2+\ve}}$. There are two residues in between, one from the pole of $D(s)$, the other from the pole of $\frac{1}{s}$. We denote them by $R_1$ and $R_2$:
\begin{align*}
R_1 &= \frac{L(1,\Sym^n(f))^{\frac{1}{n+1}}x^{1-\beta}}{(1-\beta)\beta(1+\beta)} \\
R_2 &= \frac{D(\beta)}{2}
\end{align*}
We now choose $x$ so that $x^{1-\beta} = e$. That is, choose $x = e^{\frac{\log(k)^{2+\ve}}{c_n}}$. So for $R_1$ we have
$$R_1 \ll \frac{L(1,\Sym^n(f))^{\frac{1}{n+1}}}{1-\beta}$$
For $R_2$, consider the function $(s-1)D(s)$, which is non-vanishing in the interval $(1-\frac{2c_n}{\log k},\infty)$ and positive for $s>1$. By continuity it is also positive in $(1-\frac{2c_n}{\log k},1)$. In particular $(\beta-1)D(\beta) > 0$. So $D(\beta) < 0$.
\par
We then denote 
$$I_i = \frac{1}{2\pi i}\int_{\gamma_i} \frac{D(s+\beta)x^s}{s(s+1)(s+2)}ds, i=1,2,3$$
and show they all tend to 0 as $k \to \infty$:
\begin{itemize}
\item
On $\gamma_1$, we have $\sigma = 1-\beta$ and $|t| \geq k-1$. So by convexity bound we have
\begin{align*}
D(s+\beta) &\ll \sqrt{|tk|} \\
|x^s| &= x^{1-\beta} = e \\
\frac{1}{s} &\ll \frac{1}{k}
\end{align*}
Therefore we have
$$I_1 \ll \int_{\gamma_1} \frac{\sqrt{|tk|}}{k(s+1)(s+2)}ds \ll \frac{1}{\sqrt{k}}\int_{\gamma_1}\frac{ds}{(s+1)(s+2)} \to 0 (k \to \infty)$$
\item
The convexity bound can also be used to show $I_2 \to 0 (k \to \infty)$. We skip the details here.
\item
On $\gamma_3$, same bounds hold for $D(s+\beta)$, but $\frac{1}{s}$ need to be estimated:
$$\frac{1}{s} \ll \frac{1}{\frac{c_n}{(\log k)^{2+\epsilon}} - \frac{c_n}{\log k}} \ll \log k$$
And $x^s$ is bounded by
\begin{align*}
|x^s| &= e^{\frac{(\log k)^{2+\ve}}{c_n}(\frac{c_n}{(\log k)^{2+\ve}}-\frac{c_n}{\log k})} \\
&= e^{1-(\log k)^{1+\ve}} \\
&\ll \frac{1}{e^{(\log k)^{1+\ve}}}
\end{align*}
So for $I_3$ we have
$$I_3 \ll \frac{k^2\log k}{e^{(\log k)^{1+\ve}}} \to 0 (k \to \infty)$$
\end{itemize}
Finally, by Cauchy's Theorem we have
$$I - I_1 - I_2 - I_3 = R_1 + R_2$$
and by the estimates above we have
$$1 \ll \frac{L(1,\Sym^n(f))^{\frac{1}{n+1}}}{1-\beta}$$
which gives 
$$L(1,\Sym^n(f)) \gg (1-\beta)^{n+1} \gg \frac{1}{(\log k)^{2n+2+(2n+2)\ve}}$$
where all implied constants depend on $n$ and $\ve$. This finishes the proof of Corollary 1.1.
\end{proof}

\printbibliography

@article {NT1,
    AUTHOR = {Newton, James and Thorne, Jack A.},
     TITLE = {Symmetric power functoriality for holomorphic modular forms},
   JOURNAL = {Publ. Math. Inst. Hautes \'{E}tudes Sci.},
  FJOURNAL = {Publications Math\'{e}matiques. Institut de Hautes \'{E}tudes
              Scientifiques},
    VOLUME = {134},
      YEAR = {2021},
     PAGES = {1--116},
      ISSN = {0073-8301},
   MRCLASS = {11F80 (11F70)},
  MRNUMBER = {4349240},
MRREVIEWER = {Hengfei Lu},
}

@incollection {IS,
    AUTHOR = {Iwaniec, H. and Sarnak, P.},
     TITLE = {Perspectives on the analytic theory of {$L$}-functions},
      NOTE = {GAFA 2000 (Tel Aviv, 1999)},
   JOURNAL = {Geom. Funct. Anal.},
  FJOURNAL = {Geometric and Functional Analysis},
      YEAR = {2000},
    NUMBER = {Special Volume, Part II},
     PAGES = {705--741},
      ISSN = {1016-443X},
   MRCLASS = {11M26 (11F67 11F70 11F72 11G40)},
  MRNUMBER = {1826269},
MRREVIEWER = {Henry H. Kim},
}

@book {D,
    AUTHOR = {Davenport, Harold},
     TITLE = {Multiplicative number theory},
    SERIES = {Graduate Texts in Mathematics},
    VOLUME = {74},
   EDITION = {Third},
      NOTE = {Revised and with a preface by Hugh L. Montgomery},
 PUBLISHER = {Springer-Verlag, New York},
      YEAR = {2000},
     PAGES = {xiv+177},
      ISBN = {0-387-95097-4},
   MRCLASS = {11-02 (11-01 11Mxx 11Nxx)},
  MRNUMBER = {1790423},
}

@article {HL,
    AUTHOR = {Hoffstein, Jeffrey and Lockhart, Paul},
     TITLE = {Coefficients of {M}aass forms and the {S}iegel zero},
      NOTE = {With an appendix by Dorian Goldfeld, Hoffstein and Daniel
              Lieman},
   JOURNAL = {Ann. of Math. (2)},
  FJOURNAL = {Annals of Mathematics. Second Series},
    VOLUME = {140},
      YEAR = {1994},
    NUMBER = {1},
     PAGES = {161--176},
      ISSN = {0003-486X},
   MRCLASS = {11F37 (11F12 11F30 11F66)},
  MRNUMBER = {1289494},
MRREVIEWER = {Lynne H. Walling},
}

@article {GHL,
    AUTHOR = {Dorian Goldfeld and Jeffrey Hoffstein and Daniel Lieman},
     TITLE = {Appendix: An Effective Zero-Free Region},
   JOURNAL = {Ann. of Math. (2)},
  FJOURNAL = {Annals of Mathematics. Second Series},
    VOLUME = {140},
      YEAR = {1994},
    NUMBER = {1},
     PAGES = {177--181},
      ISSN = {0003-486X},
}

@article {HR,
    AUTHOR = {Hoffstein, Jeffrey and Ramakrishnan, Dinakar},
     TITLE = {Siegel zeros and cusp forms},
   JOURNAL = {Internat. Math. Res. Notices},
  FJOURNAL = {International Mathematics Research Notices},
      YEAR = {1995},
    NUMBER = {6},
     PAGES = {279--308},
      ISSN = {1073-7928},
   MRCLASS = {11F67 (11F70)},
  MRNUMBER = {1344349},
MRREVIEWER = {David Ginzburg},
}

@article {B,
    AUTHOR = {Banks, William D.},
     TITLE = {Twisted symmetric-square {$L$}-functions and the nonexistence
              of {S}iegel zeros on {${\rm GL}(3)$}},
   JOURNAL = {Duke Math. J.},
  FJOURNAL = {Duke Mathematical Journal},
    VOLUME = {87},
      YEAR = {1997},
    NUMBER = {2},
     PAGES = {343--353},
      ISSN = {0012-7094},
   MRCLASS = {11F70 (11F66 22E55)},
  MRNUMBER = {1443531},
MRREVIEWER = {David Ginzburg},
}

@article {RW,
    AUTHOR = {Ramakrishnan, Dinakar and Wang, Song},
     TITLE = {On the exceptional zeros of {R}ankin-{S}elberg
              {$L$}-functions},
   JOURNAL = {Compositio Math.},
  FJOURNAL = {Compositio Mathematica},
    VOLUME = {135},
      YEAR = {2003},
    NUMBER = {2},
     PAGES = {211--244},
      ISSN = {0010-437X},
   MRCLASS = {11F67 (11F70 22E55)},
  MRNUMBER = {1955318},
MRREVIEWER = {Solomon Friedberg},
}

@book{IK,
    author = {H. Iwaniec and E. Kowalski},
     title = {Analytic number theory},
    series = {American Mathematical Society Colloquium Publications},
    volume = {53},
 publisher = {American Mathematical Society, Providence, RI},
      year = {2004},
      ISBN = {0-8218-3633-1},
}

@article {LW,
    AUTHOR = {Lau, Yuk-Kam and Wu, Jie},
     TITLE = {A density theorem on automorphic {$L$}-functions and some
              applications},
   JOURNAL = {Trans. Amer. Math. Soc.},
  FJOURNAL = {Transactions of the American Mathematical Society},
    VOLUME = {358},
      YEAR = {2006},
    NUMBER = {1},
     PAGES = {441--472},
      ISSN = {0002-9947},
   MRCLASS = {11F67 (11F30)},
  MRNUMBER = {2171241},
MRREVIEWER = {Emmanuel P. Royer},
}

@article {CM,
    AUTHOR = {Cogdell, J. and Michel, P.},
     TITLE = {On the complex moments of symmetric power {$L$}-functions at
              {$s=1$}},
   JOURNAL = {Int. Math. Res. Not.},
  FJOURNAL = {International Mathematics Research Notices},
      YEAR = {2004},
    NUMBER = {31},
     PAGES = {1561--1617},
      ISSN = {1073-7928},
   MRCLASS = {11F67 (11F66)},
  MRNUMBER = {2035301},
MRREVIEWER = {Emmanuel P. Royer},
}

@article {W,
    AUTHOR = {Wang, Song},
     TITLE = {On the symmetric powers of cusp forms on {${\rm GL}(2)$} of
              icosahedral type},
   JOURNAL = {Int. Math. Res. Not.},
  FJOURNAL = {International Mathematics Research Notices},
      YEAR = {2003},
    NUMBER = {44},
     PAGES = {2373--2390},
      ISSN = {1073-7928},
   MRCLASS = {11F80 (11F70 11R39 22E55)},
  MRNUMBER = {2003828},
MRREVIEWER = {Solomon Friedberg},
}

@article {NT2,
    AUTHOR = {Newton, James and Thorne, Jack A.},
     TITLE = {Symmetric power functoriality for holomorphic modular forms,
              {II}},
   JOURNAL = {Publ. Math. Inst. Hautes \'{E}tudes Sci.},
  FJOURNAL = {Publications Math\'{e}matiques. Institut de Hautes \'{E}tudes
              Scientifiques},
    VOLUME = {134},
      YEAR = {2021},
     PAGES = {117--152},
      ISSN = {0073-8301},
   MRCLASS = {11F80 (11F67 11F70)},
  MRNUMBER = {4349241},
MRREVIEWER = {Ian Kiming},
}

\end{document}